\documentclass{article}
\usepackage{epsfig}

\newtheorem{theorem}{Theorem}[section]
\newtheorem{definition}{Definition}

\begin{document}

\title{On three dimensional stellar manifolds}

\author{Sergey Nikitin \\
Department of Mathematics\\
and\\
Statistics\\
Arizona State University\\
Tempe, AZ 85287-1804}

\maketitle

{\bf Abstract} It is well known that a three dimensional (closed, connected and compact) manifold is obtained by identifying boundary faces from a stellar ball $a\star S$. The study of $S/\sim , $ two dimensional stellar sphere $S$ with $2$-simplexes identified in pairs leads us to the following conclusion: either a three dimensional manifold is homeomorphic to a sphere or to a  stellar ball $a\star S$ with its boundary $2$-simplexes identified in pairs so that $S/\sim$ is a finite number of internally flat complexes attached to a finite graph that contains at least one closed circuit.  Each of those internally flat complexes is obtained from a polygon where each side may be identified with one or more different other sides. Moreover, Euler characteristic of  $S/\sim $ is equal to one and the fundamental group of $S/\sim $ is not trivial.

\vspace{0.1cm}

\section{Introduction}

\vspace{0.1cm}

It is well known (see e.g. \cite{Massey}) that a compact $2$-dimensional surface is homeomorphic to a polygon with the edges identified in pairs.
 In this paper we show that any connected $n$-dimensional stellar manifold $M$ with a finite number of vertices possesses a similar property: $M$ is stellar equivalent to 
$$
       a\star (S/\simeq),
$$
where $a \notin S$ is a vertex, $S$ is $(n-1)$-dimensional stellar sphere. Throughout this paper it is tacitly assumed that the term "regular equivalence" (denoted as "$\simeq$") actually comprises two equivalence relations: one on the set of vertices and the other on the set of $(n-1)$-dimensional simplexes of $S.$ Those two equivalence relations satisfy the following conditions:
\begin{itemize}
\item[] No $(n-1)$-dimensional simplex of $S$ has two vertices that are equivalent to each other.
\item[] For any $(n-1)$-dimensional simplex $g$ of $S$ there might exist not more than one (different from $g$) $(n-1)$-dimensional simplex $p$ in $S$ such that$p\simeq g$ and any vertex of $g$ is equivalent to some vertex of $p.$
\end{itemize}
$a \star (S/\simeq)$ is not any more a stellar manifold in its usual sense.
We take the ball $a\star S$ and identify equivalent simplexes from $S$ but,
at the same time, we distinguish any two simplexes $a\star p$ and $a\star g$ even when $p\simeq g.$ In the sequel, 
it is convenient to use the notation $a\star (S/\simeq)$ for the resulted manifold even if it is slightly misleading. $a\star (S/\simeq)$ is called a stellar structure of $M.$\\
The complex $S/\simeq$ inherits certain properties from $a\star (S/\simeq).$

\section{Preliminaries from stellar theory}
We begin with recalling the basic definitions of stellar theory \cite{Glaser}, \cite{Lickorish}.
A stellar $n$-manifold $M$ can be identified with the sum of its $n$-dimensional simplexes ($n$-simplexes):
$$
M=\sum_{i=1}^n g_i
$$
with coefficients from ${\rm Z}_2.$ We will call $\{g_i\}_{i=1}^n$ generators of $M.$

All vertices in $M$ can be enumerated and any $n$-simplex $s$ from $M$ corresponds to the set
of its vertices
$$
s=(i_1 \; i_2 \; \dots \; i_{n+1}),
$$
where $i_1 \; i_2 \; \dots \; i_{n+1}$ are integers. 

The boundary operator $\partial $ is defined on a simplex as
$$
\partial (i_1 \; i_2 \; \dots \; i_{n+1}) = (i_1 \; i_2 \; \dots \; i_n ) + (i_1 \; i_2 \; \dots \; i_{n-1} \; i_{n+1} ) + \dots + (i_2 \; i_3 \; \dots \; i_{n+1})
$$
and linearly extended to any complex, i.e.
$$
\partial M = \sum_{i=1}^n \partial g_i.
$$
A manifold is called closed if $\partial M =0.$

If two simplexes $(i_1 \; i_2 \; \dots i_m) $ and $(j_1 \; j_2 \; \dots j_n)$ do not have common vertices then one can define their join
$$
(i_1 \; i_2 \; \dots i_m) \star (j_1 \; j_2 \; \dots j_n)
$$
as the union
$$
(i_1 \; i_2 \; \dots i_m) \cup (j_1 \; j_2 \; \dots j_n).
$$

\vspace{0.1cm}

If two complexes $K=\sum_i q_i $ and $ L = \sum_j p_j$ do not have common vertices then 
their join is defined as
$$
K\star L = \sum_{i,j} q_i \star p_j.
$$

If $A$ is a simplex in a complex $K$ then we can introduce its link:
$$
lk(A,K) = \{ B \in K \; ; \; A \star B \in K \}.
$$
The star of $A$ in $K$ is $A \star lk(A,K).$ Thus,
$$
K = A\star lk(A,K) + Q(A,K),
$$
where the complex $Q(A,K)$ is composed of all the generators of $K$ that do not contain $A.$ A complex with generators of the same dimension is called a uniform complex.

\begin{definition} ({\bf Subdivision})
Let $A$ be a simplex of a complex $K.$ Then any integer $a$ which is not a vertex of $K$ defines starring of
$$
K=A\star lk(A,K) + Q(A,K)
$$
at $a$ as
$$
 \hat{K}=a\star \partial A \star lk(A,K) + Q(A,K).
$$
This is denoted as
$$
\hat{K} =  (A \; a)K.
$$
\end{definition}

The next operation is the inverse of subdivision. It is called a stellar weld and defined as follows.

\begin{definition} ({\bf Weld})
  Consider a complex
  $$
        \hat{K}=a\star  lk(a,\hat{K}) + Q(a,\hat{K}),
  $$
with $lk(a,\hat{K}) = \partial A \star B$ where
$B$ is a subcomplex in $\hat{K},$  $A$ is a simplex and $A\notin \hat{K} .$
Then the (stellar) weld $(A\; a)^{-1} \hat{K}$ is defined as
  $$
       (A\; a)^{-1} \hat{K} =  A \star B  + Q(a,\hat{K}).
  $$
\end{definition}

A stellar move is one of the following operations: subdivision, weld, enumeration change on the set of vertices.
Two complexes $M$ and $L$  are called stellar equivalent if one is obtained from the other by a finite sequence of stellar moves.
It is denoted as $M \sim L.$ We also say that $M$ admits triangulation $L.$\\

If a complex $L$ is stellar equivalent to $(1 \; 2 \; \dots \; n+1)$ then $L$ is called a stellar 
$n$-ball. On the other hand, if $K \sim \partial (1 \; 2 \; \dots \; n+2)$ then $K$ is a stellar
$n$-sphere.

\begin{definition} ({\bf Stellar manifold})
\label{stellar_def}
Let $M$ be a complex. If, for every vertex $i$ of $M,$ the link $lk(i,M)$ is either a stellar $(n-1)$-ball or a stellar $(n-1)$-sphere, then $M$ is  a
stellar $n$-dimensional manifold ($n$-manifold).
\end{definition}

If $i$ is a vertex of $M$ then
$$
M= i\star lk(i,M) + Q(i,M).
$$
If $\partial M=0,$ then $Q(i,M)$ is a stellar manifold.\\
Indeed, consider an arbitrary vertex $j$ of $Q(i,M).$ Then
$$
lk(i,M)=j\star lk(j,lk(i,M)) + Q(j,lk(i,M))
$$
and
$$
Q(i,M)=j\star lk(j,Q(i,M)) + Q(j,Q(i,M)).
$$
Since $M$ is a stellar manifold and $\partial M = 0$
$$
i\star lk(j,lk(i,M)) + lk(j,Q(i,M)) 
$$
is a stellar sphere. Hence, it follows from \cite{Newman} that $lk(j,Q(i,M))$ is either a stellar ball or a stellar sphere.\\

The following theorem plays one of the central roles in the research on stellar manifolds.  

\begin{theorem} (Alexander \cite{Alexander})
\label{Alexander}
Let $M$ be a stellar $n$-manifold, let $J$ be a stellar $n$-ball. Suppose that
$M\cap J = \partial M \cap \partial J$ and that this intersection is a stellar $(n-1)$-ball. Then $M\cup J$ is stellar equivalent to $M.$
\end{theorem}

One can say that $M$ is obtained from $M\cup J$ by collapsing a stellar $n$-ball. If a  stellar $n$-manifold $L$ is obtained from  a stellar $n$-manifold $M$ by collapsing a finite number of simplexes (not necessary having dimension $n)$ then $M$ collapses to $L.$  If $L$ is a point then $M$ is called collapsible. \\
It was shown by Whitehead \cite{Whitehead} that a collapsible stellar manifold is a stellar ball. This fact is of special interest to us because it implies the Poincar\'e conjecture for flat stellar manifolds. 

For reader convenience we present here a definition of prism over a complex \cite{Alexandrov}.
Let $K$ be a uniform $n$-dimensional complex with vertices $i_1, \; i_2,  \; \dots \;  i_n .$ Let us take new vertices $j_1,  \; j_2 ,\; \dots \; j_n$ which correspond one-to-one to 
the vertices $i_1 , \; i_2, \; \dots \; i_n .$
The complex $P(K)$ is called a prism over $K$ if $P(K)$ consists, by definition, of all nonempty subsets of sets of the form
$$
 i_1 ,\; i_2 , \; \dots  ,\; i_k, \; j_k ,\; \dots ,\;j_r,
$$
where 
$$
 (i_1 ,\; i_2 , \; \dots  ,\; i_r) \in K.
$$
Moreover, $P(K)$ contains all the simplexes $(i_1 ,\; i_2 , \; \dots  ,\; i_r)$ from $K$ and all the simplexes $(j_1 ,\; j_2 , \; \dots  ,\; j_r)$ corresponding to them.\\
The prism $P(K)$ can be also constructed in the following way. Consider 
$$
a \star K,
$$
where $a\notin K$ is a vertex. Take its subdivision defined as
$$
L=\prod_{b \in K} ((a\; b)\; c_b) K,
$$
where $\prod_{b \in K} ((a\; b)\; c_b)$ is a superposition of subdivisions $((a\; b)\; c_b),$ where $c_b \notin K$ and $c_b \not= c_d $ for $b \not= d.$ Then
$$
L = a \star lk(a,L) + Q(a,L),
$$
where $Q(a,L)=P(K)$ is a prism over $K$ and $lk(a,L)$ is obtained from $K$ by the enumeration change $b \rightarrow c_b$ for $b\in K.$\\

In the sequel it is convenient to consider two equivalence relations: one on the set of vertices  and the other on the set of generators of a stellar manifold. Among all possible such equivalence relations we are mostly interested in those that meet certain regularity properties underlined by the following definition.
\begin{definition} ({\bf Regular equivalence})
Given a stellar manifold $M,$ a pair of equivalence relations,  one on the set of vertices and the other on the set of generators from $M,$ is called regular equivalence if it meets the following conditions:
\begin{itemize}
\item[(i)] No generator $g \in M$ has two vertices that are equivalent to each other.
\item[(ii)]For any generator $g\in M$ there might exist not more than one generator $p\in M $ such that $p \not= g$ and $p$ is equivalent to $g,\;\;g \simeq p.$ Moreover, any vertex of $g$ is either equal or equivalent to some vertex of $p.$ 
\end{itemize}
\end{definition}

Throughout the paper $a\star (S/\simeq)$ denotes the structure obtained from $a\star S$ when we identify any two equivalent generators from $S$ but distinguish the corresponding generators from $a\star S.$

\section{Stellar Manifold Structure }
We begin with proving that any connected stellar $n$-manifold with a finite number of generators can be obtained from a stellar $n$-ball $B$ by identifying in pairs the generators of $\partial B.$ For $3$-manifolds this fact was stated in \cite{SeifThrelf}.

\begin{theorem}({\bf Stellar Structure})
\label{triangulation}
A connected stellar $n$-manifold $M$ with a finite number of generators admits a triangulation (stellar structure)
$$
N = a\star (S/\simeq) ,
$$
where $a\notin S$ is a vertex, $S$ is a stellar $(n-1)$-sphere and "$\simeq $" is a regular equivalence relation. Moreover, if $M$ is closed then for any generator $g\in S$ there exists exactly one generator $p\in S$ such that $p\not= g$ and $g \simeq p.$ 
\end{theorem}
{\bf Proof.}
Let us choose an arbitrary generator $g \in M$ and an integer $a$ that is not a vertex of $M.$ Then
$$
M \sim (g \; a ) M  \mbox{  and  }
(g \; a ) M =a \star \partial g + M\setminus g,
$$
where $M\setminus g$ is defined by all the generators of $M$ excluding $g.$
We construct $N$ in a finite number of steps.
Let $N_0=(g \; a ) M.$ Suppose we constructed already $N_k$ and there exists a generator $p\in  Q(a,N_k)$ that has at least one common $(n-1)$-simplex with $lk(a, N_k).$
Without loss of generality, we can assume that 
$$
 p = (1 \; 2\; \dots \; n+1).
$$
and $(1 \; 2\; \dots \; n)$ belongs to $lk(a, N_k).$
If the vertex $(n+1)$ does not belong to  $lk(a, N_k)$ then
$$
N_{k+1} = ((a\; n+1 )\; b)^{-1}((1 \; 2 \; \dots \; n) \; b) N_k,
$$
where $b\notin N_k.$
If the vertex $(n+1)$ belongs to $lk(a, N_k)$ then after introducing a new vertex $d\notin N_k$ we take
$$
L = ((a\; d )\; b)^{-1}((1 \; 2 \; \dots \; n) \; b)(N_k\setminus p + (1\;2\;\dots \; n \;d)),
$$
where $b\notin (N_k\setminus p + (1\;2\;\dots \; n \;d )),$ and
$$
N_{k+1} =  a \star (lk(a,L)/\simeq ) + Q(a,L)
$$
endowed with the equivalence $d \simeq (n+1).$\\
By construction
$$
N_{k+1} = a \star (lk(a,N_k) +\partial p ) + Q(a,N_k) \setminus p
$$
if $(n+1)\notin lk(a, N_k).$ Otherwise,
$$
N_{k+1} = a \star ((lk(a,N_k) +\partial g)/\simeq ) + Q(a,N_k) \setminus p,
$$
where $g=(1 \; 2\; \dots \; n \; d),\;\;d\simeq (n+1).$

Since $M$ is connected and has a finite number of generators there exists a natural number $m$ such that
$$
N_m =  a \star (S/\simeq)
$$
where $S$ is a stellar $(n-1)$-sphere and "$\simeq $" is a regular equivalence relation.\\
If $M$ is closed, then $\partial N_m = 0,$ and therefore, for any generator $g \in S$ there exists exactly one generator $p\in S \setminus g$ such that $g \simeq p.$\\
Q.E.D.\\

When working with stellar manifolds we use topological terminology in the following sense. We say that a stellar $n$-manifold $M$ has a certain topological property if its standard realization in ${\rm R}^{2n+1}$ (see e.g \cite{Alexandrov}) has this property.\\
Let $\pi (M)$ denote the fundamental group of a manifold $M.$
Then Theorem \ref{triangulation} together with Seifert -- Van Kampen theorem (see e.g. \cite{Massey}) lead us to the next result.

\begin{theorem}
\label{fundamental}
If $M \sim a \star (S/\simeq )$ is a connected stellar $n$-manifold and $n>2$ then
$$
\pi (M) = \pi (S/\simeq ),
$$
where $S$ is a stellar $(n-1)$-sphere and "$\simeq $" is a regular equivalence relation.
\end{theorem}
{\bf Proof.}
By Theorem \ref{triangulation}  $M$ admits a triangulation
$$
N = a\star (S/\simeq ),
$$
where $a\notin S$ is a vertex, $S$ is a stellar $(n-1)$-sphere and "$\simeq $" is a regular equivalence relation.\\
Using the triangulation $N$ we can define open subsets $U$ and $V$ in $M$ as follows
\begin{eqnarray*}
U &=& a\star (S/\simeq )\setminus a \\
V &=& a\star (S/\simeq ) \setminus (S/\simeq ).
\end{eqnarray*}
Clearly,
$$
M = U \cup V
$$
and $V$ is homeomorphic to an interior of the $n$-ball. Thus, $V$ is simply connected and by Seifert -- Van Kampen theorem there exists an epimorphism
$$
\psi \;:\; \pi(U)\; \longrightarrow \; \pi(M)
$$
induced by inclusion $U \subset M .$ The kernel of $\psi $ is the smallest normal subgroup containing the image of the homomorphism
$$
\varphi\;:\;\pi(U\cap V) \; \longrightarrow \; \pi(U)
$$
induced by inclusion $U\cap V \subset U.$\\

The open set $U\cap V $ is homeomorphic to 
$$
0 < x_1^2+ x_2^2 +\dots +x_n^2 <1,\;\;\mbox{ where } \;\; (x_1,\; x_2 ,\; \dots ,\;x_n) \in {\rm R}^n.
$$
The fundamental subgroup of this manifold is trivial for $n>2.$ Hence, the kernel of the epimorphism $\psi $ is trivial and
$$
\pi(M) = \pi(U).
$$
On the other hand, $S/\simeq $ is a deformation retract of $U,$ and therefore,
$$
 \pi(U) = \pi (S/\simeq).
$$
Q.E.D.\\

Let $\chi (M)$ denote Euler characteristic of $M,$ i.e. 
$$
\chi(M) = \sum_{i=0}^n (-1)^i q_i,
$$
where $q_i$ is the number of $i$-simplexes in $M.$

\begin{theorem}
\label{Euler}
If $M \sim a \star (S/\simeq )$ is a closed connected stellar $n$-manifold then 
$$
      \chi(S/\simeq ) = \chi(M) + (-1)^{n+1},
$$
where $S$ is a stellar $(n-1)$-sphere and "$\simeq $" is a regular equivalence relation.
\end{theorem}
{\bf Proof.}
Since $S$ is a stellar $(n-1)$-sphere we have
$$
\chi (S) = \sum_{i=0}^{n-1} (-1)^i s_i = (-1)^{n-1}+1,
$$
where $s_i$ is the number of $i$-simplexes in $S.$
On the other hand,
$$
\chi( a \star (S/\simeq ))=\sum_{i=0}^n (-1)^i q_i,
$$
where $q_i$ is the number of $i$-simplexes in $a \star (S/\simeq ) .$ 
Let $h_i$ denote the number of $i$-simplexes in $S/\simeq .$
Taking into account that "$\simeq$" is a regular equivalence relation we obtain
$$
s_{n-1} = 2 h_{n-1},\;\;q_n = 2 h_{n-1} \;\; \mbox{ and } \;\; q_i = h_i + s_{i-1}\;\;\mbox{ for }\;\; 1\le i \le n-1
$$
and
$$
q_0 = h_0 + 1.
$$
Thus,
$$
\chi(M) = \sum_{i=0}^n (-1)^i q_i= 2 (-1)^n h_{n-1}  + h_0 + 1+ \sum_{i=1}^{n-1} (-1)^i(h_i + s_{i-1})
$$
and
$$
 \sum_{i=1}^{n-1} (-1)^i s_{i-1} = (-1)^{n-1} s_{n-1}+ (-1)^n-1 = (-1)^{n-1} 2 h_{n-1} + (-1)^n-1.
$$
Hence,
$$
\chi(M) =  \chi(S/\simeq ) + (-1)^n
$$
and the assertion is proved.\\
Q.E.D.\\

\section{Stratification of stellar manifolds}
Any stellar structure corresponds to a group defined as follows.
\begin{definition}
 Consider a closed stellar structure $a\star (S/\simeq).$ Let ${\cal G}=\{g_1,\;\dots ,g_m  \}$ be the set of generators of the stellar $n$-sphere $S.$ The regular equivalence $\simeq$ defines the permutation $p_0$ on $\cal G,$ 
$$
  p_0(g_i) = g_j \mbox{ for } g_i \simeq g_j \mbox{ and } g_i \not= g_j.
$$
We call two $(n-1)$-simplexes $\alpha $ and $\beta $ equivalent if they belong to equivalent $n$-simplexes in $S$ and each vertex of  $\alpha $ is either equal or equivalent to some vertex in $\beta .$ Let $\cal E$ denote all equivalence classes of $(n-1)$-simplexes in $S$.  Then each $\alpha \in {\cal E}$ corresponds to a permutation $p_{\alpha }$ defined as follows. If $g_i \in {\cal G}$ does not have any $(n-1)$-simplexes from $\alpha$ then
$$
    p_\alpha (g_i) = g_i.
$$
Otherwise, $g_i$ may have only one $(n-1)$-simplex from $\alpha.$ There exists only one $g_j \in  {\cal G}$ that contains the same $(n-1)$-simplex. We define 
$$
   p_\alpha (g_i) = g_j.
$$
The group of $S/\simeq$ is generated by $\{p_0,\;\{p_\alpha \}_{\alpha \in {\cal E}}\}.$
\end{definition}

In the sequel $G(S/\simeq)$ denotes the group of the stellar structure $a \star (S/\simeq).$ 
In terms of $G(S/\simeq)$ one can define the degree of a stellar structure.
\begin{definition}
 Let
$$
order(\alpha)
$$
denote the smallest integer $k>0$ such that
$$
(p_0 \cdot p_\alpha)^k = 1,
$$
Let $a\star (S/\simeq)$ be a stellar structure. Its degree $deg(S/\simeq)$ is defined as the finite string of integers
$$
deg(S/\simeq) = (deg_1(S/\simeq),\;deg_2(S/\simeq),\; \dots \; deg_m(S/\simeq))
$$
where
$$
deg_1 (S/\simeq) = \max_{\alpha \in S/\simeq} \{order(\alpha)\}
$$
and
$$
deg_j (S/\simeq) = \max_{\alpha \in S/\simeq} \{order(\alpha);\;\;order(\alpha) < deg_{j-1} (S/\simeq)\} \;\;j=2,\; \dots \; m.
$$
\end{definition}
Neither $deg(S/\simeq)$ nor $G(S/\simeq)$ are invariants of a stellar manifold.
However, these concepts lead us to the degree of a stellar manifold.
\begin{definition}
Let $M$ be a stellar manifold. Its degree
$$
deg(M)
$$
is defined as the smallest degree (in  lexicographical sense) of its stellar structures. In other words,
$$
deg(M) = \min_{(a\star S/\simeq) \sim M } \{deg( S/\simeq)\}
$$
where the minimum is taken in the sense of the lexicographical order.
\end{definition}
$deg(M)$ is an invariant of a stellar manifold. That means two stellar equivalent manifolds have the same degree. Thus, we obtain a stratification of stellar manifolds. All stellar manifolds are split into strata corresponding to different values of $deg(M).$ Stellar manifolds of the second degree $(deg(M)=(2))$ are called flat. The next theorem justifies their flatness.

\begin{theorem}
\label{flatTheorem}
 Let $M=a\star (S/\simeq )$ be an $(n+1)$-dimensional closed stellar manifold. Then $deg(S/\simeq ) = (2)$ if, and only if, each equivalence class $\alpha   $ contains not more than two different $(n-1)$-simplexes. 
\end{theorem}

{\bf Proof.}
 If  $\alpha $ contains a single $(n-1)$-simplex $i=(i_1,\dots , i_{n-1})$ then there exist two vertices $\nu$ and $\mu$ such that
\begin{eqnarray*}
         p_0(\nu \star i) &=& \mu \star i\\
         p_\alpha (\nu \star i) &=& \mu \star i
\end{eqnarray*}
 If $g$ is an $n$-simplex different from $\nu \star i$ and $ \mu \star i$ then $p_\alpha(g)=g$ and
$$
          (p_0\cdot p_\alpha)^2(g)= p^2_0(g)=g.
$$ 
On the other hand,
$$
 (p_0 \cdot p_\alpha)(\nu \star i)= \nu \star i,
$$
and
$$
 (p_0 \cdot p_\alpha)(\mu \star i)= \mu \star i.
$$
Hence, we established $(p_\alpha \cdot p_0)^2=1$ when $\alpha $ contains a single $(n-1)$-simplex. \\

If $\alpha $ contains two different $(n-1)$-simplexes $i$ and $j$ then there exist vertices $v_i$, $w_i$, $v_j$ and $w_j$ such that for the following generators of $S$ 
\begin{equation}
\label{eq_class}
v_i\star i,\; w_i \star i,\; v_j \star j,\; w_j\star j 
\end{equation}
we have
\begin{eqnarray*}
 p_\alpha (v_i\star i) &=& w_i \star i \\
 p_\alpha (v_j\star j) &=& w_j \star j \\
 p_0 (v_j\star j) &=& v_i\star i.
\end{eqnarray*}
Since $\alpha $ has only two elements we obtain
$$
p_0 (w_j\star j) = w_i\star i.
$$
If a generator $g$ is not one of the simplexes from (\ref{eq_class}) then 
$p_\alpha (g)=g,$ and $(p_0\cdot p_\alpha )^2 (g)=g.$ 
On the other hand,
$$
(p_0\cdot p_\alpha )^2 (v_i\star i)= p_0\cdot p_\alpha \cdot p_0 (w_i \star i)=p_0\cdot p_\alpha (w_j\star j) = p_0 (v_j\star j)= v_i\star i.
$$
Similar calculation show that
\begin{eqnarray*}
(p_0\cdot p_\alpha )^2 (w_j\star j) & =& w_j\star j \\
(p_0\cdot p_\alpha )^2 (v_j\star j) & =& v_j\star j \\
(p_0\cdot p_\alpha )^2 (w_i\star i) & =& w_i\star i .
\end{eqnarray*}
Hence,  $(p_\alpha \cdot p_0)^2=1$ is established when $\alpha $ contains not more than two different simplexes.\\

Now we need to show that  $(p_\alpha \cdot p_0)^2=1$ implies that the equivalence class $\alpha$ contains not more than two different $(n-1)$-simplex. We conduct the proof by reductio ad absurdum. Suppose  $\alpha$ contains at least three different  $(n-1)$-simplexes: $i,$ $j$ and $k.$ Then there exist three vertices $\mu, \;\;\nu $ and $\gamma$ such that $\mu \star i, \;\; \nu \star j$ and $\gamma \star k$ are generators from $S$ and
\begin{eqnarray*}
 p_0(\mu \star i)&=& \nu \star j\\
 p_0\cdot p_\alpha (\nu \star j)&=& \gamma \star k.
\end{eqnarray*}
Hence,
$$
 (p_\alpha \cdot p_0)^2(\mu \star i) = p_\alpha \cdot p_0 \cdot p_\alpha( \nu \star j)= p_\alpha (\gamma \star k).
$$
However,
$$
p_\alpha (\gamma \star k) \not= \mu \star i.
$$
Thus, 
$$
 (p_\alpha \cdot p_0)^2 \not= 1.
$$

Q.E.D.\\

If $(p_0 \cdot p_\alpha )^2=1$ then $(n-1)$-simplex $\alpha$ may belong to a single $n$-simplex or can be shared by two $n$-simplexes of $S/\simeq.$  We underline this fact by calling $S/\simeq$ flat at $\alpha .$

\section{Flat stellar $3$-manifolds}
Three dimensional flat stellar manifolds or stellar manifolds of the second degree admit complete classification.

\begin{theorem}({\bf Classification of flat $3$-manifolds})
\label{involution}
If $M$ is a closed connected flat stellar $3$-manifold then it is stellar equivalent to
$$
a\star (S/\simeq ), 
$$
where 
$$
S/\simeq = \left\{ \begin{array}{cc}
                         \mbox{disk } D_2 & \mbox{ for } \pi (M) = \{1\} \\
                         \mbox{projective plane } P_2 & \mbox{ for } \pi (M) = {\bf Z}_2
                  \end{array}
           \right.
$$
\end{theorem}
{\bf Proof.}
If $M$ is a flat stellar manifold  then the involution $p_0$ corresponds to a piecewise linear one-to-one involution on $S^2.$ It follows from  \cite{Lopez} that there exists a circle $S^1$ that splits $S^2$ into two parts $S_+$ and $S_-$ such that
$$
p_0:\; S^1 \to S^1
$$
and $p_0$ is a one-to-one mapping that interchanges $S_+$ and $S_-.$
By Theorem \ref{Euler}
$$
\chi (S/\simeq) =1.
$$
Hence, $S/\simeq $ is either a disk $D_2$ or a projective plane $P_2$ (see, e.g., \cite{Massey} for details).
Q.E.D.\\
 The classification of flat stellar manifolds allows to tackle certain difficult problems. In particular, the next statement proves the famous Poincare conjecture for flat stellar manifolds.
 
\begin{theorem}({\bf Poincare conjecture for flat stellar manifolds}) 
\label{Ponchik}
If $M$ is compact, closed, connected and simply connected $3$-dimensional stellar flat manifold then $M$ is a stellar $3$-sphere. 
\end{theorem}   
{\bf Proof.}
If $M$ is a closed connected flat $3$-manifold and $\pi(M)=\{1\}$ then by Theorem \ref{involution}  $S/\simeq$ is a disk.
Since  $S/\simeq $ is collapsible then so is the prism $P(S/\simeq ).$  Due to \cite{Whitehead}  the prism $P(S/\simeq)$ is a stellar ball.
$M$ is stellar equivalent to 
$$
L=\prod_{b \in S/\simeq} ((a\; b)\; c_b) (a\star (S/\simeq )),
$$
where $\prod_{b \in S/\simeq} ((a\; b)\; c_b)$ is a superposition of subdivisions $((a\; b)\; c_b),$ where $c_b \notin S/\simeq$ and $c_b \not= c_d $ for $b \not= d.$
$$
 L=  a \star \partial P(S/\simeq ) + P(S/\simeq ),
$$
$L$ is a union of two balls with identified boundaries. Hence \cite{Glaser}, $L$ is a stellar $3$-sphere and so is $M.$
Q.E.D.\\

This simple proof is not valid for stellar manifolds that are not flat.

\section{ Stellar $3$-manifolds}
Now we consider three dimensional manifolds that are not flat. Examples of such non-flat manifolds are provided by lens spaces.\\
  A lens shell $\ell (q,p)$ for coprime integers 
$$
q > p \ge 1
$$
can be defined as follows. Consider the function of a complex variable $z\in {\rm C},$
$$
p_0(z)=\left\{ \begin{array}{cc}
               \frac{e^{2\cdot  \pi \frac{p}{q} i }}{\bar z },\;\;&\mbox{ for } z\cdot \bar z \le 1,\\
               \frac{e^{-2\cdot  \pi \frac{p}{q} i}}{\bar z },\;\;&\mbox{ for } z\cdot \bar z  > 1,
               \end{array}
	\right.
$$
Now define the equivalence on $\{ {\rm C}, \infty \} $ as
$$
0 \sim \infty,\;\;\mbox{ and } z_1 \sim z_2 \;\;\mbox{as long as } z_1 = p_0^m(z_2)
$$
for some integer $m,$ where $ p_0^2(z)$ denotes $ p_0( p_0(z))$ and $p_0^m(z) = p_0(p_0^{m-1}(z)).$

$\{ {\rm C}, \infty \} $ is identified with the $2$-dimensional sphere $S$ (e.g. with stereo-graphic projection). The topological space
$$
S/\sim
$$
is called lens shell $\ell (q,p).$\\
Notice that in the literature the {\it lens space } is defined as  $a\star \ell (q,p).$ \\
It is not difficult to show that 
$$
\pi(\ell (q,p) ) = {\rm Z_q}
$$
The points 
$$
e^{2\cdot m \pi \frac{p}{q} i }\;\;m=0,\;1,\;\dots q-1
$$
are all mapped into a single vertex of $S/\sim .$ Consider $1$-simplex $\alpha $ in $S/\sim $ such that 
$$
p_\alpha (z) = \frac{1}{\bar z} .
$$
It is easy to see that
$$
(p_\alpha \cdot p_0)^q = 1
$$
In other words, for $q>2$ lens shell $\ell(q,p)$ is not flat and its degree is equal to $q.$ \\
We call an edge $\alpha \in S/\sim$ {\it collapsible } if there exists a $2$-simplex $F\in S$ such that
$$
p_\alpha (p_0 ( F )) = F.
$$
Consider the subgroup $G_2(S/\sim)$ of group $G(S/\sim)$ generated by
$$
\{p_\alpha :\;order(\alpha )= 2, \; \alpha \in (\partial P)/\sim , \;\alpha  \mbox{ is not collapsible } \}.
$$
Then the set of $2$-simplexes from $S$ is partitioned into the orbits of $G_2(S/\sim ).$ All those orbits are identified in pairs in order to obtain $S/\sim .$ For the sake of brevity we call each pair of identified orbits {\it internally flat complex.} \\
We call $1$-simplex {\it non-flat } if its order is higher than $2.$ Then the boundary of an internally flat complex consists of non-flat $1$-simplexes and  collapsible  $1$-simplexes. If the boundary contains a collapsible $1$-simplex then the internally flat complex is collapsible.\\
Taking in to account Theorem \ref{triangulation} we conclude that each internally flat complex is obtained from a polygon where each boundary side can be identified with one or more other different boundary sides.\\

\begin{theorem}
A three dimensional manifold is either homeomorphic to a sphere  or to $a\star (S/\sim)$ such that $S/\sim$ is the sum of a finite number of internally flat complexes attached to a finite graph that contains at least one closed circuit. Moreover, Euler characteristic $\chi (S/\sim)=1.$
\end{theorem}
{\bf Proof.}
The statement of this theorem is already established for flat manifolds (Theorem \ref{involution}). Our goal is to show that it is valid for for non-flat manifolds.\\
 If $M$ is not flat then it is stellar equivalent to  $a\star (S/\sim)$ with faces identified in pairs and we can consider the graph $\Gamma $ that consists of all edges $\alpha \in S/\sim $ such that $order(\alpha ) >2. $
 If $deg(M)>2$ then  $\Gamma $ is not empty. Moreover, $\Gamma $ has at least one closed circuit. Indeed, if $\Gamma $ was a set of trees then one could consider it as a part of one single spanning tree. Outside of this tree $S/\sim $ is internally flat. That means $S/\sim$ consists of a finite number of internally flat complexes attached to that tree. Each internally flat complex is  collapsible. Since $\Gamma $  is a set of trees then $S/\sim$ neither contains shell spaces nor $2$-manifolds with non-trivial genus. On the other hand, it is known (see Theorem \ref{Euler} ) that $\chi (S/\sim) =1.$ That means $S/\sim$ is collapsible, and therefore,  $M$ is homeomorphic to three dimensional sphere (see \cite{Glaser}, \cite{Whitehead} for details).  That contradicts to the assumption that $deg(M )>2.$ Thus,
 $$
 deg(M)>2
 $$
 implies the existence of at least one closed circuit in $\Gamma .$\\
Q.E.D.\\
The graph $\Gamma $ formed by non-flat $1$-simplexes from $S/\sim$ plays important role in understanding properties of $S/\sim .$ 
In particular, if $\Gamma $ is stellar equivalent to a circle then $S/\sim $ is a lens shell.\\
\begin{theorem}
\label{BigPonchik}
A stellar $3$-manifold is either stellar equivalent to $3$-sphere or its fundamental group is not trivial.
\end{theorem}
{\bf Proof.} The statement of this theorem was already proved for flat manifolds (see Theorem \ref{Ponchik}).\\
 Consider a non flat manifold $M, \;\; deg(M)>2.$ $M$ is homeomorphic to $a\star (S/\sim ).$ Without loss of generality we can assume that $S/\sim$ does not have collapsible complexes. The graph $\Gamma $ formed by non flat edges of $S/\sim $ has at least one closed circuit. Take a minimal spanning tree for $\Gamma $ (the minimal tree in $S/\sim $ that contains all vertices from $\Gamma $ and the maximal possible number of arcs from $\Gamma $) and contract it to a point, say $x_0.$ Then $\Gamma $ is homotopic to a finite sum of circles $ \bigvee_{j=1}^n S^1_j $ (wedge of circles) having the only common point $x_0.$ We conduct the proof by mathematical induction with respect to the number of circles.\\
 Since $S/\sim$ is a $CW$-complex the homotopy of $\Gamma $ can be extended to the homotopy of $S/\sim $ (see \cite{Shubert}). Hence, if the number of circles is equal to one, then $S/\sim $ is homotopic to a lens shell. This homotopy corresponds to an epimorphism of the fundamental group for  $S/\sim $ onto the fundamental group of the lens shell. Hence, the fundamental group of $S/\sim $ is not trivial.\\
 Suppose it is true that the fundamental group of $S/\sim $ is not trivial when $\Gamma $ is homotopic to the sum of $n$ circles attached to each other at the only point $x_0.$  Consider $S/\sim $ such that  $\Gamma $ is homotopic to the sum of $n+1$ circles $ \bigvee_{j=1}^{n+1} S^1_j $ attached to each other at $x_0.$ If the fundamental group was trivial then $S^1_{n+1}$ would be homotopic to $x_0$ or to the sum of circles from $ \bigvee_{j=1}^n S^1_j.$
 Moreover, this homotopy can be extended to the homotopy of  $S/\sim .$ That yields an epimorphism of the fundamental group $\pi(S/\sim )$ onto the fundamental group of the homotopy of  $S/\sim $ with $\Gamma $ homotopic to  the sum of $n$ circles. However, it contradicts to the assumption that the fundamental group is not trivial for those $S/\sim $ that have $\Gamma $ homotopic to  the sum of $n$ circles.\\
Q.E.D.\\

\bibliographystyle{ams-plain}

\begin{thebibliography}{10}
\bibitem{Alexander} J.W. Alexander, \textit{The combinatorial theory of complexes,} Ann, of Math. {\bf 31} (1930) 292-320.
\bibitem{Alexandrov} P.S. Alexandrov, \textit{Combinatorial Topology}, Volumes I, II and III, Dover Publications, Inc., 1998.
\bibitem{Glaser} L. C. Glaser, \textit{Geometrical Combinatorial Topology}, Volume I, Van Nostrand Reinhold Company, N.Y., 1970.
\bibitem{Lickorish} W. B. R. Lickorish,  \textit{Simplical moves on complexes and manifolds}, Volume 2: Proceedings of the Kirbyfest, Geometry and Topology Monographs(1999), 299-320.
\bibitem{Lopez} S. L\`opez de Medrano , \textit{ Involutions on Manifolds }, Springer-Verlag, N.Y., Berlin, 1971.
\bibitem{Massey} W.S. Massey,  \textit{Algebraic topology: an introduction,} Graduate Texts in Math. vol.56, Springer, 1967.
\bibitem{Newman} M.H.A. Newman, \textit{On the foundations of combinatorial Analysis Situs,} Proc. Royal Acad. Amsterdam, {\bf 29} (1926) 610-641.
\bibitem{SeifThrelf}H.Seifert, W.Threlfall, {Lehrbuch der Topologie}, Chelsea Publishing Co., New York, 1945.
\bibitem{Shubert} H. Shubert, {\em Topology, } Allan and Bacon, Inc., 1968, Boston.
\bibitem{Whitehead}J. H. C. Whitehead, \textit{Mathematical works}, Pergamon, 1962.
\end{thebibliography}

\end{document}